\documentclass[11pt]{amsart}
\newtheorem{theo}{THEOREM}[section]

\newtheorem{cor}[theo]{Corollary}
\newtheorem{prop}[theo]{Proposition}
\newtheorem{dfntn}[theo]{Definition}
\newtheorem{rem}[theo]{Remark}
\newtheorem{claim}[theo]{Claim}

\newtheorem{ex}[theo]{Example}
\newtheorem{facts}[theo]{FACTS}
\newcommand{\brref}[1]{(\ref{#1})}

\newcommand{\Pin}[1]{{\mathbb P}^{#1}}
\newcommand{\restrict}[2]{{#1}_{\mid _{#2}}}

\newcommand{\Projcal}[1]{\mathbb{ P}({\mathcal #1})}
\newcommand{\calo}{{\mathcal O}}

\newcommand{\chiS}{\chi(\calo_S)}
\newcommand{\oof}[2]{{\mathcal O}_{#1}({#2})}
\newcommand{\oofp}[2]{{\mathcal O}_{\mathbb{ P}^{#1}}({#2})}

\newcommand{\lra}{\longrightarrow}

\newcommand{\scrollcal}[1]{(\Projcal{#1},\tautcal{#1})}
\newcommand {\xel} {(X, L)}

\newcommand{\tautcal}[1]{{\mathcal O}_{\mathbb{P}({\mathcal#1})}(1)}
\newcommand{\tautcalof}[2]{{\mathcal O}_{\mathbb{P}({\mathcal#1})}(#2)}

\newenvironment{rem*}{\begin{rem}\em}{\end{rem}}
\newenvironment{ex*}{\begin{ex}\em}{\end{ex}}
\newenvironment{claim*}{\begin{claim}\em}{\end{claim}}
\newenvironment{facts*}{\begin{facts}\em}{\end{facts}}

\title{The dimension of the Hilbert scheme of special threefolds\\
}
\author{Gian Mario Besana}
\address{Gian Mario Besana\\
C.T.I. DePaul University\\
243 S. Wabash \\
Chicago IL 60604\\USA}
\email{gbesana@cs.depaul.edu}

\author{Maria Lucia Fania}
\address{ Maria Lucia Fania\\ Dipartimento di Matematica
\\Universit\`{a} degli Studi di L'Aquila\\
Via Vetoio Loc. Coppito\\67100 L'Aquila\\Italy}
\email{fania@univaq.it}

\thanks{Partially supported by MIUR of the Italian Government in the
framework
of the National Research Project (Cofin 2002) {\it Geometria sulle
Variet\`a Algebriche}.\\
The material in this paper is, in part, based upon work supported by
the National Science Foundation (NSF) under Grant No. 0125068.  Any
opinions, findings and conclusions or recommendations expressed in this
material are those of the authors and do not necessarily reflect the
views of the NSF}

\begin{document}

\begin{abstract}

The Hilbert scheme of
3-folds in $\Pin{n},$ $n\ge 6,$ that are scrolls over
$\Pin{2}$ or over a smooth quadric surface $\mathbf{Q} \subset \Pin{3}$
or
that are quadric or cubic fibrations over $\Pin{1}$ is studied.
All known such threefolds of degree $7 \le d \le 11$ are shown to
correspond to smooth points of an irreducible component of their Hilbert
scheme, whose dimension is computed.

\end{abstract}

\maketitle
\section{introduction}

The classification of complex projective manifolds of low degree,
extending
the classical works of Weil, \cite{We}, and Swinnerton Dyer, \cite{SD},
has
been conducted in the recent past as a three-step process. Maximal lists
of possible manifolds are first compiled, according to the admissible
values
of their numerical invariants. The second step deals with establishing
the actual existence of manifolds in the lists, looking for effective
constructions of explicit examples. Finally, the Hilbert scheme of
existing manifolds with given Hilbert polynomial is investigated.

The first step of the program has seen over the years the successful
work
of several authors.  Manifolds of degree up to eight have been
classified
by Ionescu,
\cite{Io1}, \cite{Io2}, \cite{Io3},
Okonek, \cite{ok2}, \cite{ok8}, Abo, Decker and Sasakura, \cite{ADS}.
Other authors have considered the special case of manifolds of fixed
low codimension. Beltrametti, Schneider and Sommese, \cite{BSS1},
\cite{BSS2},
Ottaviani, \cite{ot1}, Braun, Ottaviani, Schneider and Schreyer,
\cite{boss}, classified  $3$-folds in
$\Pin{5}$ of degree up to twelve, and recently  Bertolini, \cite{ber},
classified $3$-folds of degree twelve in $ \Pin{6}$.
The classification of manifolds of dimension $n\ge 3$, regardless of the
codimension, of degree nine and ten was accomplished
by   Fania and Livorni, \cite{fa-li9}, \cite{fa-li10},  while Besana and
Biancofiore, \cite{be-bi}, considered the case of degree eleven.

Many of the authors cited above also dealt with the issue of the actual
existence of the classified manifolds, although open questions still
remain in this arena. Further contributions to this second step of the
program can also be found in works of Okonek, \cite{ok3}, and
Biancofiore
and Fania, \cite{bi-fa}.

As to the Hilbert schemes, general results are known only for
fixed low codimension. In codimension two Ellingsrud, \cite{e},
considered arithmetically Cohen Macaulay varieties regardless of
their degree, Ionescu dealt with such cases among the manifolds in
\cite{Io2}, while the Hilbert scheme of a special class of
$3$-folds in $\Pin{5}$ was studied by Fania and Mezzetti,
\cite{fa-me}. General results in codimension two are also due to
M.C. Chang, \cite{chang1}, \cite{chang2}.

In codimension $3$, Kleppe and Mir\'o-Roig, \cite{kl-mr}, computed
the dimension of the open subset of the Hilbert scheme of
arithmetically Gorenstein closed subschemes, independently of
their degree, while Kleppe, Migliore, Mir\'o-Roig, Nagel and
Peterson, \cite{kmmnp}, dealt with good determinantal subschemes.

In this paper, four classes of $3$-folds in  $\Pin{n}$ with $n \ge 6$
are
considered.
  Precisely we deal with  $3$-folds which are scrolls over $\Pin{2}$ or
  over a smooth quadric surface $\mathbf{Q}$ and
$3$-folds which are quadric or cubic fibrations over $\Pin{1}.$
Their geometric special structure  is exploited in order to find mild
conditions
guaranteeing that these manifolds are unobstructed and their Hilbert
scheme has an irreducible component which is smooth at the point
corresponding
to the given variety. Explicit formulas to compute the dimension of such
components are given.

All manifolds in the above classes, of degree $7 \le d \le 11,$ and
known to
exist, are shown to be unobstructed. The dimension of an irreducible
component of their Hilbert scheme is computed. In some cases, a dense
open subset of such an irreducible component of the Hilbert scheme is
shown to be the locus of good determinantal subschemes with the given
Hilbert polynomial (see Section \ref{gooddet} for definitions).

The paper is structured as follows: in Section \ref{notation} we
collect all the necessary notation and background material,
including the fundamental theorem of Grothendieck on the existence
of the Hilbert scheme; in Section \ref{scrolls} unobstruction
results for scrolls over $\Pin{2}$ and $\mathbf{Q}$ are presented;
Section \ref{fibrations} is devoted to the same results for
quadric and cubic  fibrations over $\Pin{1};$ Section
\ref{gooddet} explores the connections with results on
unobstructed good determinantal subschemes presented in
\cite{kmmnp}.

\section{Notation and Preliminaries}
\label{notation}
In this section notation is fixed and definitions and
results which will be used throughout the paper are
recalled.

  Let $X$  be a complex projective manifold of
dimension $3,$  $3$-fold for short, and let $L$ be a very
ample line bundle on $X$. Projective properties of  $\xel$ are always
referred to the embedding $X \subset \Pin{n}$ given by the complete
linear system associated
with
$L.$  Therefore our $3$-folds are always  linearly normal.
  We denote by  $S$ and $C$, respectively, a smooth
  surface and curve,  obtained as transverse
intersection of $X$ with respectively  $1$ and  $2$  general elements
of $|L|$.
  For any coherent sheaf ${\mathcal F}$ on $X$,  $h^i(X,{\mathcal F})$
is the complex dimension of
$H^i(X,{\mathcal F}).$ When the ambient variety is understood, we often
write $H^i({\mathcal F})$
and $h^i({\mathcal F})$ respectively for $H^i(X,{\mathcal F})$ and
$h^i(X,{\mathcal F}).$ The following
notation will be used throughout this work.
  \begin{enumerate}
\item[ ]$\chi(L) = \sum(-1)^ i h^i(L)$, the Euler characteristic of $L$;
\item[ ]  $\restrict{L}{Y}$ the restriction of $L$ to a subvariety $Y;$
\item[ ]  $K_X$ the
canonical bundle of $X.$ When the context is clear,
$X$ may be dropped;
\item[ ] $q(S) = h^1 ({\mathcal
O}_{S})$, the irregularity of $S$;
\item[ ] $p_g
(S) = h^0 (K_S)$, the geometric genus of $S$;
\item[
] $c_i = c_i(X)$,  the $i^{th}$ Chern class
of $X$;
\item[ ] $d = \deg{X} = L^3$, the degree of $X$ in the embedding
given by $L$;
\item[] $g = g(X),$ the sectional genus of $\xel$ defined by
$2g-2=(K+2L)L^2;$
\item[] $\mathbf{Q},$ a smooth quadric hypersurface embedded in
$\Pin{3};$
\item[] $\lceil x \rceil,$ the ceiling of  a real number $x,$ i.e. the
smallest integer greater than,  or equal to, $x.$
\end{enumerate}

  Cartier
divisors, their associated line bundles and the invertible sheaves
of their holomorphic sections are used with no distinction. Mostly
additive notation is used for their group. Multiplicative
notation (juxtaposition) will be used for intersection  of cycles and
Chern classes.

%
%
\begin{dfntn}
\label{specialvar} A pair $(X, L)$, where
$L$ is an ample line bundle on a $3$-fold $X,$ is
a {\it scroll}, or a {\it hyperquadric fibration}, or  a {\it Del
Pezzo fibration} over a normal variety $Y$ if there exist an ample line
bundle $M$  on $Y$ and a surjective morphism  $\varphi: X \to Y$ with
connected fibers
such that $K_X + (4 - \dim Y) L = \varphi^*(M)$  or,
respectively,
  $K_X + (3 - \dim Y) L = \varphi ^*(M)$ or
  $K_X + (2 - \dim Y) L = \varphi^*(M).$
\end{dfntn}
%
%
\begin{rem*}
\label{resonscronsup}
Let $\xel$ be a scroll over a surface $Y.$ Then, see \cite{BESO} 14.1.3,
$X \cong \Projcal{E} $, where ${\mathcal E}= \varphi_{*}(L)$ and
$L$ is   the tautological  line bundle on $\Projcal{E}.$ It is
known, see for example \cite{BESO} Section 11.1, that $S$ is the blow
up of $Y$
at $c_2(\mathcal{E})$ points, $\chi({\mathcal O}_{Y}) =
\chi({\mathcal O}_{S})$ and $d = c_1^2(\mathcal{E})-c_2(\mathcal{E}).$
\end{rem*}

%
%
%

\begin{rem*}
  \label{resonhqf}
Let $(X,L)$ be  a $3$-dimensional manifold which is either a
hyperquadric fibration or a Del Pezzo fibration of fiber degree $3$
over ${\mathbb P}^1,$
$\varphi:X \to {\mathbb P}^1,$  which is embedded  by $|L|$  in
${\mathbb P}^{n}$.
Let $F\in|\varphi^*({\mathcal O}_{{\mathbb P}^1}(1))|$ be a fiber of
$\varphi$,
then $L^2F = \alpha = 2,3,$ respectively. As $L$ is very ample,
\cite{Io1},\cite{fug2} and  \cite{ds1} show that there is a rank $4$
vector bundle over ${\mathbb P}^1$,
$\mathcal{E}=\varphi_{*}(L)={\mathcal O}_{{\mathbb P}^1}(a_4)\oplus
\dots \oplus {\mathcal O}_{{\mathbb P}^1}(a_1)$ and an embedding
$\iota: X \to \Projcal{E},$ such that
$L = \iota^*(\tautcal{E}).$ Note that $X$ is a divisor on
$\Projcal{E}$ and  $X\in |\tautcalof{E}{\alpha}   + \rho^* \oofp{1}{b}|$
  where $\rho : \Projcal{E} \to \Pin{1}$ is
the projection map. Let $ e = \deg
{\mathcal E},$   then
\begin{xalignat}{3} \label{hqfs}
\begin{cases}
d = \alpha e + b\\
2 g - 2 = \alpha( e + b -2) + (\alpha -2) d
\end{cases} &  &\text{or }&  & \begin{cases}
b=\frac {d - d \alpha + 2 g - 2 + 2 \alpha}{\alpha - 1}\\
e = - \frac {2(g - 1 + \alpha - d \alpha + d)}{\alpha(\alpha - 1)}
\end{cases} &
\end{xalignat}
\end{rem*}
\vskip .2 in
\subsection{The Hilbert Scheme}
The existence of the Hilbert scheme for closed subschemes of $\Pin{n}$
with given Hilbert polynomial was established by Grothendieck,
\cite{groth}. The following formulation of his basic result is due to
Sommese, \cite{SO1}.

\begin{prop}[\cite{groth}, \cite{SO1}]
\label{basic fact} Let $Z$ be a smooth connected projective variety.
Let $X$ be a connected submanifold of $Z$ with $H^1(X,N)=0$ where
$N$ is the normal bundle of $X$.
Then there exist irreducible projective varieties ${\mathcal Y}$ and
${\mathcal H}$ with the following properties:
\begin{itemize}
\label{properties}
  \item[(i)]
${\mathcal Y} \subset {\mathcal H}\times Z$ and the map
$p:{\mathcal Y}\to {\mathcal H}$ induced
by the product projection is a flat surjection,
  \item[(ii)]  there is a smooth point $x\in {\mathcal H}$
with $p$ of maximal rank in a neighborhood of $p^{-1}(x)$,
  \item[(iii)]  $q$ identifies  $p^{-1}(x)$ with $X$ where
$q:{\mathcal Y}\to Z$ is the map induced
by the product projection, and
\item[(iv)] $H^0(N)$ is naturally identified with $T_{{\mathcal H},x}$
where
$T_{{\mathcal H},x}$ is the Zariski tangent space of ${\mathcal H}$ at
$x$.
\end{itemize}
\end{prop}

%
%

\section{Hilbert scheme of 3-dimensional scrolls over $\Pin{2}$  or
$\mathbf{Q}$}
\label{scrolls}
Let $\xel = \scrollcal{E}$ be  a smooth 3-fold of sectional
genus $g$ and degree $d$ which is a scroll over ${\mathbb P}^2,$ as in
Definition \ref{specialvar} and Remark \ref{resonscronsup}. Let $X$ be
embedded
by $|L|$ in $\Pin{n}.$
The following proposition shows that, under mild conditions on the
embedding
and on the splitting type of $\mathcal{E},$ $X$ is unobstructed.

\vspace{1.5mm}

\begin{prop}
\label{Hilbertscheme of scrollover P2}
Let $\xel = \scrollcal{E}$ be a 3-dimensional scroll over $\Pin{2}$ of
degree $d$
and sectional genus $g.$ Let $X$ be embedded by $|L|$ in $\Pin{n}.$ Let
$c_1({\mathcal E})={\mathcal O}_{\Pin{2}}(e_1)$, $e_2=c_2({\mathcal
E}).$ Assume:
\begin{itemize}
\item[i)] $H^1(X,L)=0;$
\item[ii)] there exists a line $\ell \subset \Pin{2}$ such that
$\restrict{\mathcal{E}}{\ell} = \oofp{1}{a} \oplus \oofp{1}{e_1 -a}$
where $\frac{e_1}{2} \le a \le \frac{e_1}{2} + 1.$
\end{itemize}
Then the  Hilbert scheme of $X$ has an irreducible component, ${\mathcal
H}$,  which is smooth at the point
representing $X$ and
$$\dim{\mathcal
H}= (d+2) (n-3) +\frac{3e_1}{2}(n+1)-\frac{e_1^2}{2}(n-5) -4$$
\end{prop}
\begin{proof}
Let $N$ denote the normal bundle of $X$ in ${\mathbb P}^n$.
  The statement will follow from Proposition \ref{basic fact} by showing
that $H^1(X,N)=0$ and conducting an explicit computation of $h^0(X,N).$
Let
\begin{eqnarray}
\label{eulersequscrollsuP2}
0\lra {\mathcal O}_{X} \lra {\mathcal O}_{X}(1)^{\oplus (n+1)}
\lra T_{{{\mathbb P}^n}|{X}} \lra  0
\end{eqnarray}
be the Euler sequence on ${\mathbb P}^n$
restricted to $X$. As $(X,L)$ is a scroll over ${\mathbb P}^2$,
\begin{eqnarray}
\label{quadratino}
H^{i}(X,{\mathcal O}_{X})= H^{i}({\mathbb P}^2,{\mathcal O}_{{\mathbb
P}^2})= 0,\quad \text{ for}\quad i\ge 1.
\end{eqnarray}
Let $S\in |L|$ and $C\in |L_S|$ be a general surface and curve section.
 From their structure sequences, noting that $S$ is rational, it follows
that
\begin{eqnarray}
\label{doppioquadratino}
H^2(X,L)= H^3(X,L)= 0.
\end{eqnarray}
Because $H^1(X,L)=0$ by assumption, from \brref{quadratino},
\brref{doppioquadratino} and the cohomology sequence associated to
\brref{eulersequscrollsuP2} it follows that  $H^i(X,T_{{ {\mathbb
P}^n}|{X}})=0$
for $i\ge 1$. Therefore the exact sequence
\begin{eqnarray}
\label{tangentsequ}
0\lra T_{X} \lra T_{{ {\mathbb P}^n}|{X}}
\lra N \lra 0
\end{eqnarray}
gives
\begin{eqnarray}
\label{cohomnormal}
H^{i}(X,N) = H^{i+1}(X,T_{X}) \qquad {\text {for} \quad  i\ge 1.}
\end{eqnarray}
In particular $H^3(X,N) = 0$ for dimension reasons.
To compute
$H^{j}(X,T_{X}),j = 2,3,$
let $\varphi:{\mathbb P}({\mathcal E})\lra {\mathbb P}^2$ be the scroll
map,
and consider the relative cotangent bundle sequence:

\begin{eqnarray}
\label{relativctgbdl}
0\to \varphi^{*}({\Omega}^1_{{\mathbb P}^2})\to {\Omega}^1_{X}
\to {\Omega}^1_{X|{{\mathbb P}^2}} \lra 0.
\end{eqnarray}

 From \brref{relativctgbdl} and the Whitney sum one obtains
$$ c_1({\Omega}^1_{X})= c_1(\varphi^{*}({\Omega}^1_{{\mathbb
P}^2}))+c_1({\Omega}^1_{X|{{\mathbb P}^2}})$$ and thus
$${\Omega}^1_{X|{{\mathbb P}^2}}=K_X+\varphi^{*}({\mathcal
O}_{\Pin{2}}(3)).$$

The adjunction theoretic characterization of the scroll then gives
$${\Omega}^1_{X|{{\mathbb P}^2}}=K_X+\varphi^{*}({\mathcal
O}_{\Pin{2}}(3))=
-2L+\varphi^{*}({\mathcal O}_{\Pin{2}}(e_1)) $$

that, combined with the dual of \brref{relativctgbdl}, gives

\begin{eqnarray}
\label{relative tgbdl}
0 \to 2L-\varphi^{*}({\mathcal O}_{\Pin{2}}(e_1)) \to T_X
\to \varphi^{*}(T_{\Pin{2}}) \to 0.
\end{eqnarray}

As the cohomology of $\varphi^{*}(T_{\Pin{2}})$ is easily computed, we
devote our attention to the cohomology of
$2L-\varphi^{*}({\mathcal O}_{\Pin{2}}(e_1)).$
Noticing that  $R^{i}\varphi_{*}(2L)=0$ for $i \ge 1$
(see \cite{H}, pg 253), projection formula
and Leray's spectral
sequence give
$$H^i(X, 2L-\varphi^{*}({\mathcal O}_{\Pin{2}}(e_1)))\cong
H^i({\Pin{2}}, S^2{\mathcal E}\otimes {\mathcal O}_{\Pin{2}}(-e_1)).$$
Therefore
\begin{equation}
\label{h3sym2E}
H^3(X, 2L-\varphi^{*}({\mathcal O}_{\Pin{2}}(e_1))) = 0
\end{equation} for dimension reasons.

Let $\ell$ be a line  in $\Pin{2}$ as in assumption ii).

It is
$$({\mathcal E}\otimes {\mathcal E})_{|\ell}=\oofp{1}{2a}\oplus
\oofp{1}{2(e_1-a)}\oplus \oofp{1}{e_1} \oplus \oofp{1}{e_1}.$$
As $det({\mathcal E})=\oofp{2}{e_1}$ and ${\mathcal E}\otimes {\mathcal
E}=
det({\mathcal E})\oplus S^2{\mathcal E}$ it follows that
\begin{equation}
\label{restrictionofs2E}
S^2{\mathcal E}_{|\ell}=
\oofp{1}{2a}\oplus \oofp{1}{2(e_1-a)} \oplus \oofp{1}{e_1}.
\end{equation}

Tensoring the structure sequence of $\ell$ on $\Pin{2}$ with $
S^2{\mathcal E}\otimes \oofp{2}{t}$  gives
\begin{equation}
\label{2ndsymE}
0 \to S^2{\mathcal E}(t-1) \lra S^2{\mathcal E}(t) \to
\restrict{(S^2{\mathcal E}\otimes \oofp{2}{t})}{\ell}  \to 0
\end{equation}
where
\begin{equation}
\label{twisteds2E}
\restrict{(S^2{\mathcal E}\otimes \oofp{2}{t})}{\ell} =
\oofp{1}{2a+t}\oplus \oofp{1}{2(e_1-a)+t} \oplus \oofp{1}{e_1+t}.
\end{equation}
Because $a \ge \frac{e_1}{2}$ one can see that the minimum of the
integers $2a, 2(e_1-a), e_1$ is $2(e_1-a)$.

The cohomology sequence associated with  \brref{2ndsymE} then gives

  $$H^2(S^2{\mathcal E}(t-1))=H^2(S^2{\mathcal E}(t)) \quad  \forall \,
t \ge 2(a-e_1)-1.$$
Thus

$$H^2(S^2{\mathcal E}(2(a-e_1)-2))=H^2(S^2{\mathcal
E}(2(a-e_1)-1))=\dots
= H^2(S^2{\mathcal E}(s)) = \dots.$$

Therefore Serre's vanishing theorem gives
$$H^2(S^2{\mathcal E}(s))=0, \quad \text{for all } s \ge 2(a-e_1)-2.$$
In particular $H^2(S^2{\mathcal E}(-e_1))=0$ because of our assumption
$a \le \frac{e_1}{2}+1,$
and thus, from  the cohomology
sequence associated to \brref{relative tgbdl},\brref{h3sym2E}, and
dimension reasons, it follows that
$H^2(X, T_X)= H^2(\Pin{2}, f^*T_{\Pin{2}})$ and $H^3(X, T_X) = 0.$
On the other hand, by Leray spectral sequence,
$$H^2(\Pin{2}, \varphi^*T_{\Pin{2}}) = H^2(\Pin{2}, T_{\Pin{2}})=0.$$
Hence $H^i(X,T_{X})=0, i= 2,3$ and thus, by \brref{cohomnormal},
$$H^i(X,N)=0 \quad i = 1,2,3.$$
According to Proposition \ref{basic fact}, there exists an irreducible
component $\mathcal{H}$ of the Hilbert scheme of $X \subset \Pin{n}.$
The dimension of ${\mathcal H}$, by Proposition \ref{basic fact},
  (iv), will be given by $h^0(X,N)= \chi(N).$
The Hirzebruch-Riemann-Roch theorem gives
\begin{eqnarray}
\label{chiN}
  \chi(N) &=& \frac{1}{6}(n_1^3-3n_1n_2+3n_3)+
\frac{1}{4}c_1(n_1^2-2n_2)+\frac{1}{12}(c_1^2+c_2)n_1 \\
& & +(n-3)\chi({\mathcal O}_{X})\nonumber
\end{eqnarray}
where $n_i=c_i(N),$ and $ c_i=c_i(X).$

Chern classes
of $N$  can be obtained from
\brref{tangentsequ}:
\begin{xalignat}{1}
\label{valueniscrollP}
n_1=&K+(n+1)L;  \\
n_2=&\frac{1}{2}n(n+1)L^2+(n+1)LK+K^2-c_2; \notag\\
n_3=&\frac{1}{6}(n-1)n(n+1)L^3+\frac{1}{2}n(n+1)KL^2+
(n+1)K^2L-(n+1)c_2L \notag \\
&-2c_2K+K^3-c_3. \notag
\end{xalignat}

The numerical invariants of $X$ can be easily computed:
\begin{xalignat}{2}
  KL^2 &= -2d+e_1^2-3e_1;&  K^2L &= 4d-3e_1^2+6e_1+9;\notag \\
  c_2L &= 3e_1+3;& K^3 &= -8d+6e_1^2-54;\notag \\
  -Kc_2 &= 24;&  c_3 &= 6. \nonumber
\end{xalignat}

  Plugging these in \brref{valueniscrollP} and the results in
\brref{chiN} one gets
$$\chi(N)=-10+\frac{3}{2}e_1-3d+\frac{5}{2}e_1^2+2n+\frac{3}{2}ne_1+dn-
\frac{1}{2}ne_1^2.$$
\end{proof}

%
%

Let $\xel=\scrollcal{E}$ be  a smooth 3-fold of sectional
genus $g$ and degree $d$
which is a scroll over a smooth quadric surface $
\mathbf{Q}\subset {\mathbb P}^3,$
  as in Definition \ref{specialvar} and Remark \ref{resonscronsup}. Let
$X$ be
embedded by $|L|$ in $\Pin{n}.$
The following proposition shows that, under mild conditions
on the embedding and on the splitting type of $\mathcal{E},$ $X$
is unobstructed.

\vspace{1.5mm}

\begin{prop}
\label{Hilbertscheme of scrollover Q2}
Let $\xel = \scrollcal{E}$ be a 3-dimensional scroll over $\mathbf{Q}$
of degree $d$
and sectional genus $g.$ Let $X$ be embedded by $|L|$ in $\Pin{n}$ and
let
$c_1({\mathcal E})= \oof{\mathbf{Q}}{e_{11},e_{12}}$, $c_2({\mathcal
E})= e_2.$  Assume:
\begin{itemize}
\item[i)] $H^1(X,L) = 0$;
\item[ii)]there exists a line
$\ell_1\in |{\mathcal O}(1,0)|$ and a line $\ell_2\in |{\mathcal
O}(0,1)|$
such that:
$\restrict{\mathcal
E}{\ell_{1}}=\oofp{1}{\lceil\frac{e_{12}}{2}\rceil}\oplus
\oofp{1}{e_{12}-\lceil\frac{e_{12}}{2}\rceil},$\\
$\restrict{\mathcal
E}{\ell_2}=\oofp{1}{\lceil\frac{e_{11}}{2}\rceil}\oplus
\oofp{1}{e_{11}-\lceil\frac{e_{11}}{2}\rceil}.$
\end{itemize}
Then the  Hilbert scheme of $X$ has an irreducible component, ${\mathcal
H}$,  which is smooth at the point
representing $X$ and
$$\dim{\mathcal
H}=(d+2)(n-3)+(e_{11}+e_{12})(n+1) - e_{11}e_{12}(n-5)-2.$$
\end{prop}
\begin{proof}

The proof proceeds exactly as in Proposition \ref{Hilbertscheme of
scrollover P2}. Letting $N$ denote the normal bundle of $X$ in
${\mathbb P}^n,$ the computation of $H^i(N)$ relies this time on
$H^i(\mathbf{Q}, S^2{\mathcal E}\otimes {\mathcal
O}_{\mathbf{Q}}(-e_{11},-e_{12})).$

Let $\ell_1$ and $\ell_2$ be lines in $|\oof{\mathbf{Q}}{1,0}|$ and
$|\oof{\mathbf{Q}}{0,1}|$,
respectively.

We can assume that $\restrict{\mathcal E}{\ell_1}=\oof{\ell_1}{a}\oplus
\oof{\ell_1}{e_{12}-a}$  and
$\restrict{\mathcal E}{\ell_2}=\oof{\ell_2}{b}\oplus
\oof{\ell_2}{e_{11}-b}$ with
$a\ge \frac{e_{12}}{2}$ and $b\ge \frac{e_{11}}{2}.$

Tensoring the structure sequences of $\ell_1$ and $\ell_2$ on
$\mathbf{Q}$ by
$S^2{\mathcal E}(t-e, t-1)$ and $S^2{\mathcal E}(t-e-1, t)$
respectively, where $e=e_{11}-e_{12},$ we get:

\begin{gather}
\label{2ndsymE1twisted}
0\to S^2{\mathcal E}(t-e-1,t-1) \to S^2{\mathcal E}(t-e,t-1) \to \\
{\mathcal O}_{\ell_{1}}(2a+t-1)\oplus
  {\mathcal O}_{\ell_{1}}(e_{12}+t-1)
  \oplus{\mathcal O}_{\ell_{1}}(2(e_{12}-a)+t-1) \to 0,\nonumber\\
\nonumber \\
\label{2ndsymE2twisted}
0 \to S^2{\mathcal E}(t-e-1,t-1) \to S^2{\mathcal E}(t-e-1,t)  \to \\
{\mathcal O}_{\ell_{2}}(2b+t-e-1)\oplus{\mathcal
O}_{\ell_{2}}(e_{11}+t-e-1)\oplus{\mathcal
O}_{\ell_{2}}(2(e_{11}-b)+t-e-1)
\to 0.\nonumber
\end{gather}

Assume
  \begin{equation}
  \label{ell1vanishing}
  2(e_{12}-a)+t-1\ge -1,
  \end{equation}
   that is $t\ge 2(a-e_{12})=2(a-e_{11}+e),$ and assume
\begin{equation}
\label{ell2vanishing}
2(e_{11}-b)+t-e-1 \ge -1,
\end{equation}
  that is $t\ge 2(b-e_{11})+e.$
Then the cohomology sequences associated with \brref{2ndsymE1twisted}
and \brref{2ndsymE2twisted} give

  \begin{equation}
  \label{catenah2combinata}
  H^2(S^2{\mathcal E}(t-e,t-1))=H^2(S^2{\mathcal E}(t-e-1,t-1))=
  H^2(S^2{\mathcal E}(t-e-1,t)).
  \end{equation}

Tensoring now the structure sequence of $\ell_1$ with $S^2{\mathcal
E}(t-e,t)$ one gets:
\begin{eqnarray}
\label{2ndsymE1twisted1}
0 \lra S^2{\mathcal E}(t-e-1,t) \lra S^2{\mathcal E}(t-e,t) \lra
{\mathcal O}_{\ell_{1}}(2a+t)\oplus\\
\oplus
  {\mathcal O}_{\ell_{1}}(e_{12}+t)\oplus{\mathcal
O}_{\ell_{1}}(2(e_{12}-a)+t) \lra 0\nonumber
\end{eqnarray}
 From \brref{ell1vanishing} it follows:

\begin{equation}
\label{catenah2ell1}
H^2(S^2{\mathcal E}(t-e,t))=H^2(S^2{\mathcal E}(t-e-1,t)).
\end{equation}
Thus, under assumptions \brref{ell1vanishing} and
\brref{ell2vanishing}, from \brref{catenah2combinata} and
\brref{catenah2ell1}, it follows that
$$H^2(S^2{\mathcal E}(t-e,t))=H^2(S^2{\mathcal E}(t-e-1,t-1)),$$
thus by Serre's vanishing theorem, under the same assumptions,
$$H^2(S^2{\mathcal E}(t-e-1,t-1)) =0.$$

To obtain the desired vanishing $H^2(S^2{\mathcal
E}(-e_{11},-e_{12}))=0$ the following conditions need to be satisfied:

\begin{itemize}
\label{conditions}
\item[(a)] $2(a-e_{12})-e-1 \le -e_{11};$
\item[(b)]  $2(b-e_{11})+e-1 \le -e_{12}.$
\end{itemize}

Condition (a) along with the fact that $a\ge \frac{e_{12}}{2}$ gives
$a=\lceil\frac{e_{12}}{2}\rceil.$
Analogously, condition (b) gives $b=\lceil\frac{e_{11}}{2}\rceil.$

Similarly to the proof of Proposition \ref{Hilbertscheme of scrollover
P2} we get
$H^i(X,T_X) = 0$ for $i=2,3,$ thus by \brref{cohomnormal}
$$H^i(X,N) = 0 \quad i= 1,2,3,$$ and therefore, by Proposition
\ref{basic fact}, we have established the existence of an irreducible
component $\mathcal{H}$ of the Hilbert scheme of $X \subset \Pin{n}.$

The dimension of ${\mathcal H},$ as in Proposition \ref{Hilbertscheme
of scrollover P2}, is obtained via Hirzebruch-Riemann-Roch theorem,
using \brref{tangentsequ} to compute the Chern classes of $N.$

Let us now consider the numerical invariants of $X.$ Let $H_i =
\varphi^*(\ell_{i})$ and  let $F$ denote a fiber of $\varphi$. One can
easily obtain the following relations in the cohomology ring of
$X={\mathbb P}({\mathcal E}):$
\begin{xalignat}{4}
\label{intersectionincohomringX}
L^3& = d; & L^2H_1& = e_{12}; & L^2H_2& = e_{11};&  H_1^3& = H_2^3 =
0;\\
L^2F& = 1; & LH_1H_2& = 1; & H_1F& = H_2F = 0 ;
   & LH_1^2&=LH_2^2=0. \nonumber
\end{xalignat}
Using \brref{intersectionincohomringX} we  get:
\begin{xalignat}{2}
\label{valueinvarscrollQ}
KL^2&=-2d+2(e_{11}e_{12}-e_{11}-e_{12});&  c_2L &=
2(e_{11}+e_{12})+4;\nonumber\\
K^2L&=4d+4(e_{11}+e_{12})-6e_{11}e_{12}+8;& K^3 &=
-8d+12e_{11}e_{12}-48; \nonumber\\
  -Kc_2& = 24;&& \nonumber
\end{xalignat}

and thus, recalling that $c_3=8,$
$$h^0(X,N)=\chi(N)=n(d+2+e_{11}+e_{12}-e_{11}e_{12})-8-
3d+e_{11}+e_{12}+5e_{11}e_{12}.$$

\end{proof}
%
%
%
As a Corollary to Proposition \ref{Hilbertscheme of scrollover P2} we
show
that all the known
$3$-folds scrolls over $\Pin{2}$ of
degree $7 \le d \le 12$ are unobstructed and we compute the dimension of
the irreducible component of the Hilbert scheme to which they belong.
Note that no such varieties exist for $d\le 6$.

\begin{cor}
\label{corscrollp2} Let $\xel = \scrollcal{E}$ be a 3-dimensional
scroll over $\Pin{2}.$ Let $X$ be embedded by $|L|$ in $\Pin{n},$ with
degree $d$ and sectional genus $g$
as in the table below.
Then the  Hilbert scheme of $X \subset \Pin{n}$ has an irreducible
component, ${\mathcal
H}$,  which is smooth at the point
representing $X$ and of dimension as in the rightmost column of the
table.
\vskip .2in
\begin{tabular}{|c|c|c|c|c|c|c|}
\hline
$d$& $g$& $n$&$c_1(\mathcal{E})$&$c_2(\mathcal{E})$&  Reference &
$\dim{\mathcal{H}}$\\
\hline
\hline
$7$&$3$&$6 $&$4$&$9$&\cite{Io2} Prop. 1.3&$57$\\
\hline
$8$&$3$&$7$&$4$&$8$&\cite{Io2} Prop. 1.3&$68$\\
\hline
$9$&$3$&$8$&$4$&$7$&\cite{Io2} Prop. 1.3,&$81$\\
& & & & &\cite{fa-li9} Prop. 3.1 and Remark 3.2&\\
\hline
$10$&$3$&$9$&$4$&$6$&\cite{Io2} Prop. 1.3,&$96$\\
& & & & &\cite{fa-li10} Prop. 3.4&\\
\hline
$10$&$6$&$6$&$5$&$15$&\cite{fa-li10} Remark 5.3&$72$\\
\hline
$12$&$3$&$11$&$4$&$4$&\cite{Io1} Prop. 4.7&$132$\\
\hline
\end{tabular}
\end{cor}
\begin{proof}
We will show that for all the cases in the above table, the hypothesis
of Proposition \ref{Hilbertscheme of scrollover P2} are satisfied.
The structure sequences of a general surface section $S$ and curve
section $C$ tensored with $L$ and $\restrict{L}{S}$ respectively, give
\begin{equation}
h^0(L) - h^1(L)= 3 + d - g.
\end{equation}
On the other hand, $h^0(L) = n+1$ and thus
\begin{equation}
\label{h1Lscrollsup2}
h^1(L) = n-2-d+g.
\end{equation}
A simple check gives $h^1(L)= 0$ for all the cases in the above table.

To establish the existence of a line $\ell$ as in hypothesis ii) of
Proposition \ref{Hilbertscheme of scrollover P2} we first consider the
cases in the above table with $c_1(\mathcal{E}) =4.$ As $\mathcal{E}$
is ample, the generic splitting type of $\mathcal{E}$ is then either
$\oofp{1}{2}\oplus \oofp{1}{2}$ or $\oofp{1}{3}\oplus \oofp{1}{1}.$ In
both cases, a generic line $\ell$ satisfies hypothesis ii).

In the case with $c_1(\mathcal{E}) = 5$ the possible splitting
types are $\oofp{1}{3}\oplus \oofp{1}{2}$ or $\oofp{1}{4}\oplus
\oofp{1}{1}.$ If the generic splitting type is $\oofp{1}{3}\oplus
\oofp{1}{2}$, then hypothesis ii) is satisfied for a generic line
$\ell.$ If the generic splitting type is $\oofp{1}{4}\oplus
\oofp{1}{1},$ then $\mathcal{E}$ must be uniform. If not, there
should exist a line on which $\mathcal{E}$ jumps, i.e. on which it
splits as $\oofp{1}{a}\oplus \oofp{1}{b},$ $a \ge b,$ with $(4,1)\le (a,b),$ lexicographically, see \cite[p.29]{OSS}; this is impossible
as the only other admissible splitting type for $\mathcal{E}$ is
$\oofp{1}{3}\oplus \oofp{1}{2},$ being $\mathcal{E}$ ample, and $(3,2) \le (4,1),$
lexicographically.

Uniform 2-bundles on $\Pin{2}$ either split or are of the form
$T_{\Pin{2}}(t)$, see \cite[Theorem 2.2.2, p. 211]{OSS}. Both
cases are ruled out as $c_1(\mathcal{E}) = 5$ and
$c_2(\mathcal{E}) =15.$

By Proposition \ref{Hilbertscheme of scrollover P2} there exists an
irreducible component $\mathcal{H}$ of the Hilbert scheme of $X \subset
\Pin{n}$ whose dimension can now be easily computed.
\end{proof}
\begin{rem*}
The $3$-dimensional scrolls over $\Pin{2}$ of degree $11$ whose
invariants are the following:
\begin{center}
\begin{tabular}{|c|c|c|c|c|c|c|}
\hline
$d$& $g$& $n$&$c_1(\mathcal{E})$&$c_2(\mathcal{E})$&  Reference &
$\dim{\mathcal{H}}$\\
\hline
\hline
$11$&$3$&$10$&$4$&$5$&\cite{be-bi} Prop. 4.2.2&$113$\\
\hline
$11$&$6$&$7$&$5$&$14$&\cite{be-bi} Prop. 5.2.3&$83$\\
\hline
\end{tabular}
\end{center}
also satisfy the hypothesis of Proposition \ref{Hilbertscheme of
scrollover P2}. Unfortunately the existence of these scrolls is not
known as yet.
\end{rem*}
%
%
%

As a Corollary to Proposition \ref{Hilbertscheme of scrollover Q2} we
show
that all the known
$3$-folds scrolls over $\mathbf{Q}$
of degree $8 \le d \le 11$ are unobstructed and we compute the dimension
of the irreducible component of the Hilbert scheme to which they belong.
Note that no such varieties exist for $d\le 7$.

\begin{cor} Let $\xel = \scrollcal{E}$ be a 3-dimensional scroll over
$\mathbf{Q}.$ Let $X$ be  embedded by $|L|$ in $\Pin{n},$ with degree
$d$ and sectional genus $g$
as in the table below.
Then the  Hilbert scheme of $X \subset \Pin{n}$ has an irreducible
component, ${\mathcal
H}$,  which is smooth at the point
representing $X$ and of dimension as in the rightmost column of the
table.
\vskip .2in
\begin{tabular}{|c|c|c|c|c|c|c|}
\hline
$d$& $g$& $n$&$c_1(\mathcal{E})$&$c_2(\mathcal{E})$&  Reference &
$\dim{\mathcal{H}}$\\
\hline
\hline
$8$&$4$&$6$&$\oof{\mathbf{Q}}{3,3}$&$10$&\cite{Io2} Prop. 1.6&$61$\\
\hline
$9$&$4$&$7$&$\oof{\mathbf{Q}}{3,3}$&$9$&\cite{fa-li10} Remark 7.5&$72$\\
\hline
$10$&$4$&$8$&$\oof{\mathbf{Q}}{3,3}$&$8$&\cite{fa-li10} Remark
3.5&$85$\\
\hline
$11$&$4$&$9$&$\oof{\mathbf{Q}}{3,3}$&$7$&\cite{be-bi} Remark 4.2.4
&$100$\\
\hline
\end{tabular}
\end{cor}
\begin{proof}
We will show that for all the cases in the above table, the hypothesis
of Proposition \ref{Hilbertscheme of scrollover Q2} are satisfied.
To verify hypothesis i) of Proposition \ref{Hilbertscheme of scrollover
Q2} one proceeds exactly as in the first part of the proof of Corollary
\ref{corscrollp2}.

As $c_1(\mathcal{E}) = \oof{\mathbf{Q}}{3,3},$ the splitting type of
$\mathcal{E}$ on any line of both rulings is, $\oofp{1}{2}\oplus
\oofp{1}{1},$ and thus hypothesis ii) is satisfied.

By Proposition \ref{Hilbertscheme of scrollover Q2} there exists an
irreducible component $\mathcal{H}$ of the Hilbert scheme of $X \subset
\Pin{n}$ whose dimension can now be easily computed.
\end{proof}

%
%

\section{Hilbert scheme of 3-dimensional fibrations over
$\Pin{1}$ with low fiber degree}
\label{fibrations}
In this section we deal with $3$-folds that are fibrations over
$\Pin{1}$ with
low fiber degree. As in the case of scrolls over surfaces we will see
that under mild conditions these $3$-folds are also unobstructed.
The dimension of the irreducible component of the Hilbert scheme to
which they belong is computed.

Further notation is introduced here below.
\subsection{Notation}
\label{fibrationsonP1}
Let $(X,L)$ be  a $3$-dimensional manifold which is either a
hyperquadric fibration or a Del Pezzo fibration of fiber degre $3$ over
${\mathbb P}^1,$ as in Definiiton \ref{specialvar}. Let $X$ be embedded
by $|L|$ in $\Pin{n}.$
Let $F\in|\varphi^*({\mathcal O}_{{\mathbb P}^1}(1))|$ be a fiber of
$\varphi$,
then $L^2F = \alpha = 2,3,$ respectively. Note that
$\mathcal{E}=\varphi_{*}(L)$
  is a rank $4$ vector bundle over ${\mathbb P}^1$,
$\mathcal{E}={\mathcal O}_{{\mathbb P}^1}(a_4)\oplus
\dots \oplus {\mathcal O}_{{\mathbb P}^1}(a_1).$ We can arrange the
$a_i$ so that
  $a_{4}\ge \dots \ge  a_1$.  Then  $X$ is embedded in a
${\mathbb P}^{3}$-bundle $\Projcal{E}$ over
${\mathbb P}^{1}$, $\iota: X \to \Projcal{E}$ such that $L = \iota
^*(H)$ where
$H=\tautcal{E}$ is
the tautological line bundle of $\mathcal{E}$ and
$X\in|\alpha H+\rho^*{\mathcal O}_{{\mathbb P}^1}(b)|$, for some
integer $b$, where
$\rho:\Projcal{E} \to {\mathbb P}^1$ is the projection map.
In what follows the notation is as above.

\vspace{1.5mm}

\begin{prop}
\label{Hilbertscheme of fibroverP1}
Let $(X,L)$ be  a $3$-dimensional manifold which is either a
hyperquadric
fibration or a Del Pezzo fibration over ${\mathbb P}^1.$ Let
$\mathcal{E},$ $\alpha = 2,3$, $a_i, i=1,\dots, 4,$ $b,$ $H,$ $\rho$,
$\iota,$ and $\varphi$ be as in {\rm \ref{fibrationsonP1}}.
Let $d$  be the  degree of $X$ and let
$g$ be its sectional genus. If
\begin{itemize}
\item[i)] $H^1(X,L)= 0,$
\item[ii)] $-\alpha a_1-1\le b,$
\end{itemize}
then the  Hilbert scheme of $X$ has an irreducible component, ${\mathcal
H}$,  which is smooth at the point
representing $X$ and
\begin{equation}
\label{dimfibroverP1}
\dim{\mathcal H} =
\begin{cases}
{d (n-4) + g (14 - n) + 8 + 3n }, & \text{if $\alpha = 2;$} \\
{\frac{2d}{3}(n-14)+\frac{g}{3}(44-n)+\frac{10}{3}(10+n)}, & \text{if
$\alpha = 3.$} \\
\end{cases}
\end{equation}
\end{prop}
\begin{proof}
Let $N$ denote the normal bundle of $X$ in ${\mathbb P}^{n}$.
  In order to use Proposition \ref{basic fact} we need to  show that
$H^1(X,N)=0.$  Noticing that, if $\alpha =2,3,$ $R^i\rho_*((1-\alpha)H)
=0, i\ge 0,$  the structure sequence of $X$ in $\Projcal{E}$ tensored
with $H$ gives $H^i(\Projcal{E}, H) = H^i \xel,\, i\ge 0.$ On the other
hand $H^i(\Projcal{E}, H) = H^i (\Pin{1}, E) = 0,$ for $i\ge 2,$ and
thus, recalling that $H^1(X,L) = 0$ by assumption,  $H^i(X,L) = 0, i\ge
1.$

Reasoning as in the proof of
Proposition \ref{Hilbertscheme of scrollover P2} we get that
\begin{eqnarray}
\label{cohomnormaloverP1}
h^{i}(X,N) = h^{i+1}(X,T_{X}) \qquad {\text {for} \quad  i\ge 1.}
\end{eqnarray}
Hence $H^3(X,N) = 0$ for dimension reasons.
As $X\subset \Projcal{E}$ and
$X\in|\alpha H+\rho^*{\mathcal O}_{{\mathbb P}^1}(b)|$, we
have
the following exact sequences:
\begin{eqnarray}
\label{structuresequinPE}
0\lra {\mathcal O}_{\Projcal{E}}(-X) \lra {\mathcal O}_{\Projcal{E}}
\lra {\mathcal O}_{X} \lra 0,
\end{eqnarray}
\begin{eqnarray}
\label{tangentsequinPE}
0\lra T_{X} \lra T_{{ \Projcal{E}}|{X}}
\lra (\alpha H+\rho^*{\mathcal O}_{{\mathbb P}^1}(b))_{|X} \lra 0.
\end{eqnarray}

  Note that $(\alpha H+\rho^*{\mathcal O}_{{\mathbb P}^1}(b))_{|X}
= \alpha L+\varphi^*{\mathcal O}_{{\mathbb P}^1}(b))=\alpha L+bF,$ where
  $F\in |\varphi^*{\mathcal O}_{{\mathbb P}^1}(1)|$ is a fiber of
$\varphi$.

To compute
$H^2(X,T_{X})$ is enough to compute  $H^i(X, \alpha L+bF)$ and
$H^i(X, T_{{ \Projcal{E}}|{X}}).$

  Tensoring  sequence \brref{structuresequinPE}
with $\alpha H+\rho^*{\mathcal O}_{{\mathbb P}^1}(b)$  we get
\begin{eqnarray}
\label{structuresequinPEtwisted}
0\lra {\mathcal O}_{\Projcal{E}} \lra \alpha H+\rho^*{\mathcal
O}_{{\mathbb P}^1}(b)
\lra \alpha L+\varphi^*{\mathcal O}_{{\mathbb P}^1}(b) \lra 0.
\end{eqnarray}
As $${H^i(\Projcal{E},{\mathcal O}_{\Projcal{E}})=
h^i({\mathbb P}^1,{\mathcal O}_{{\mathbb P}^1})=0} \text{ for $i\ge
1$}$$
it follows that
$$H^i(X,\alpha L+\varphi^*{\mathcal O}_{{\mathbb P}^1}(b))
=H^i(\Projcal{E},\alpha H+\rho^*{\mathcal O}_{{\mathbb P}^1}(b)) \text{
for $i\ge 1$}.$$
By Leray's spectral sequence $$H^i(\Projcal{E},\alpha H+\rho^*{\mathcal
O}_{{\mathbb P}^1}(b))=
H^i({\mathbb P}^1,S^{\alpha}(\mathcal{E})\otimes {\mathcal O}_{{\mathbb
P}^1}(b)),$$
hence $$H^i(X,\alpha L+\varphi^*{\mathcal O}_{{\mathbb P}^1}(b))=
H^i({\mathbb P}^1,S^{\alpha}(\mathcal{E})\otimes {\mathcal O}_{{\mathbb
P}^1}(b))=0\text{ for $i\ge 2$}.$$

Because the ${a_i}'s$ are sorted in increasing order, the smallest
degree of the line bundles of the decomposition of
$S^{\alpha}(\mathcal{E}) \otimes \oof{\Pin{1}}{b}$ is $\alpha a_1 + b.$
 From our assumption $-\alpha a_1-1\le b$ it follows
$$H^1(X,\alpha L+\varphi^*{\mathcal O}_{{\mathbb P}^1}(b))=
H^1({\mathbb P}^1,S^{\alpha}(\mathcal{E})\otimes {\mathcal O}_{{\mathbb
P}^1}(b))=0.$$

Tensoring sequence \brref{structuresequinPE} with $T_{\Projcal{E}}$
  we get
\begin{eqnarray}
\label{structuresequinPEtensorT}
0\to T_{\Projcal{E}}(-\alpha H+\rho^*{\mathcal O}_{{\mathbb P}^1}(-b))
\to
  T_{\Projcal{E}}
\to  T_{\Projcal{E}|X} \to 0.
\end{eqnarray}
  From  sequence \brref{tangentsequinPE} it follows that  $H^i(X,
T_{X})=H^i(T_{{\Projcal{E}}|{X}}), i \ge 2.$ Vanishing of these
cohomology groups will follow from $H^{i+1}(T_{\Projcal{E}}(-\alpha
H+\rho^*{\mathcal O}_{{\mathbb P}^1}(-b))=0,$ $H^i(T_{\Projcal{E}})=0,$
  $i\ge 2,$ and  \brref{structuresequinPEtensorT}.

In order to compute such cohomology groups we consider the
following exact sequences associated to $\rho:\Projcal{E} \to {\mathbb
P}^1:$

\begin{eqnarray}
\label{relativetgt1}
0\to T_{\Projcal{E}|{\mathbb P}^1}
\to T_{\Projcal{E}}
\to  \rho^*T_{{\mathbb P}1} \to 0,
\end{eqnarray}
\begin{eqnarray}
\label{relativetgt2}
0\to {\mathcal O}_{\Projcal{E}}
\to \rho^*E^*\otimes {\mathcal O}_{\Projcal{E}}(1)
\to  T_{\Projcal{E}|{\mathbb P}^1} \to 0.
\end{eqnarray}
Tensoring  sequences \brref{relativetgt1} and \brref{relativetgt2}
with $-\alpha H+\rho^*{\mathcal O}_{{\mathbb P}^1}(-b)$
we get, respectively,

\begin{eqnarray}
\label{relativetgt1twisted}
0\to
T_{\Projcal{E}|{\mathbb P}^1}(-\alpha H+\rho^*{\mathcal O}_{{\mathbb
P}^1}(-b))
\to \\ T_{\Projcal{E}}(-\alpha H+\rho^*{\mathcal O}_{{\mathbb P}^1}(-b))
\to
-\alpha H+\rho^*{\mathcal O}_{{\mathbb P}^1}(2-b) \to 0, \nonumber
\end{eqnarray}
\begin{eqnarray}
\label{relativetgt2twisted}
0\to -\alpha H+\rho^*{\mathcal O}_{{\mathbb P}^1}(-b)
\to \rho^*(E^*(-b))\otimes {\mathcal O}_{\Projcal{E}}(1-\alpha)
\to  \\
T_{\Projcal{E}|{\mathbb P}^1}(-\alpha H+\rho^*{\mathcal O}_{{\mathbb
P}^1}(-b))\to 0. \nonumber
\end{eqnarray}
Taking direct images of
\brref{relativetgt1twisted} via $\rho,$ and
noticing that, if $\alpha = 2,3,$ it is
$$R^i\rho_{*}(-\alpha H+\rho^*{\mathcal O}_{{\mathbb P}^1}(-b))= 0,
\quad \text{for }i\ge 0, $$
it follows that
$$R^i\rho_{*}(T_{\Projcal{E}|{\mathbb P}^1}(-\alpha H+\rho^*{\mathcal
O}_{{\mathbb P}^1}(-b)))=
R^i\rho_{*}(T_{\Projcal{E}}(-\alpha H+\rho^*{\mathcal O}_{{\mathbb
P}^1}(-b))), \quad \text{for $i\ge 0$}.$$

Similarly, taking direct images via $\rho$ of sequence
\brref{relativetgt2twisted} and noticing that
$$R^i\rho_{*}(-\alpha H+\rho^*{\mathcal O}_{{\mathbb P}^1}(-b))=
R^i\rho_{*}( \rho^*(E^*(-b))\otimes {\mathcal O}_{\Projcal{E}}(-2))=0,
\quad \text {for $i\ge 0$},$$
it follows that
$$R^i\rho_{*}(T_{\Projcal{E}|{\mathbb P}^1}
(-\alpha H+\rho^*{\mathcal O}_{{\mathbb P}^1}(-b)))=0,
  \quad \text {for $i\ge 0$}.$$
Hence $H^i(\Projcal{E}, T_{\Projcal{E}|{\mathbb P}^1}
(-\alpha H+\rho^*{\mathcal O}_{{\mathbb P}^1}(-b))=
H^i({\mathbb P}^1, \rho_{*}(T_{\Projcal{E}|{\mathbb P}^1}
(-\alpha H+\rho^*{\mathcal O}_{{\mathbb P}^1}(-b)))=0,$ for all $i\ge
0.$

We now turn our attention to $H^i(\Projcal{E}, T_{\Projcal{E}}).$
\\First notice that
$H^i(\Projcal{E}, \rho^*{\mathcal O}_{{\mathbb P}^1}(2)))=
H^i({\mathbb P}^1, {\mathcal O}_{{\mathbb P}^1}(2))=0,$ for  $i\ge 1.$
Taking direct images via $\rho$ of sequence
\brref{relativetgt2} and
noticing that
$$R^i\rho_{*}({\mathcal O}_{\Projcal{E}})=
R^i\rho_{*}(\rho^*E^*\otimes {\mathcal O}_{\Projcal{E}}(1))= 0,
\quad \text {for $i\ge 1,$}$$
it follows that
\begin{eqnarray}
\label{relativetgt2directimage}
0\to {\mathcal O}_{{\mathbb P}^1}\to E^* \otimes E
\to \rho_{*}(T_{\Projcal{E}|{\mathbb P}^1}) \to 0
\end{eqnarray}
and
$$R^i\rho_{*}(T_{\Projcal{E}|{\mathbb P}^1})= 0,
  \quad \text {for $i\ge 1$}.$$

Thus by Leray's spectral sequence
$H^i(\Projcal{E},T_{\Projcal{E}|{\mathbb P}^1})=
H^i({\mathbb P}^1,\rho_{*}(T_{\Projcal{E}|{\mathbb P}^1}))$ for  $i\ge
0.$
Moreover, as
$H^i({\mathbb P}^1,{\mathcal O}_{{\mathbb P}^1})=0$ for  $i\ge 1$ and
$H^i({\mathbb P}^1,E^* \otimes E)=0$ for  $i\ge 2,$
it follows that
$H^i({\mathbb P}^1,\rho_{*}(T_{\Projcal{E}|\Pin{1}}))= 0,$ for $i \ge
2.$
Using the cohomology sequence associated to \brref{relativetgt1}
  we get that
$H^i(\Projcal{E},T_{\Projcal{E}})=0, i\ge 2.$
Hence $H^i(X,T_{X})=0, i\ge 2,$ and thus by \brref{cohomnormaloverP1}
we get $H^1(X,N)=0, i\ge 1.$

The dimension of ${\mathcal H}$, by (\brref{basic fact},
  (iv)), is now obtained via Hirzebruch-Riemann-Roch, as in the proof of
Proposition \ref{Hilbertscheme of scrollover P2}.

One can compute easily the numerical invariants of $X$ and get:
\begin{align}
  KL^2&= 2g-2-2d; \nonumber\\
  K^2L &= \alpha d(4-\alpha)+4g(\alpha-4)-4\alpha+16; \nonumber\\
  c_2L &= \alpha d(2-\alpha)+2g(2\alpha-3)+2\alpha^2-2\alpha+6;
\nonumber \\
  -Kc_2&=24; \nonumber \\
K^3&=2d(16-24\alpha+9\alpha^2-\alpha^3)+6g(16-8\alpha+\alpha^2)-
6\alpha^2+48\alpha-96; \nonumber \\
c_3 &=
2d(\alpha^3-3\alpha^2+3\alpha-1)-6g(\alpha^2-2\alpha+1)-
4\alpha^3+10\alpha^2-6\alpha+6. \nonumber
\end{align}

Thus
\begin{eqnarray*}
h^0(N)=\chi(N)& =&
\frac{d}{6}(-\alpha^3+5\alpha^2+4n-n\alpha^2+3n\alpha+40-38\alpha)+\\
& &\frac{g}{6}(94 - 14n-11\alpha+3\alpha^2 + 4n\alpha)+\\
& &\frac{1}{6}(38\alpha+\alpha^3+20n+19\alpha^2+n\alpha^2-3n\alpha-112).
\end{eqnarray*}
\end{proof}
%
%

As a Corollary to Proposition \ref{Hilbertscheme of fibroverP1} we show
that
all the known
$3$-folds which are either a hyperquadric
fibration or a Del Pezzo fibration over ${\mathbb P}^1$
  of
degree $7 \le d \le 12$ are unobstructed and we compute the dimension of
the irreducible component of the Hilbert scheme to which they belong.
Note that no such varieties exist for $d\le 6$.

\begin{cor}
\label{corquadricfibrations} Let $\xel$ be as in {\rm
\ref{fibrationsonP1}}, with $\alpha = 2,$ i.e. a hyperquadric fibration
over $\Pin{1}.$ Let $X$ be embedded by $|L|$ in $\Pin{n}$ and let the
numerical invariants of $X \subset \Pin{n}$ be as in the table below.
Then the  Hilbert scheme of $X$ has an irreducible component, ${\mathcal
H}$,  which is smooth at the point
representing $X$ and of dimension as in the rightmost column of the
table.
\begin{center}
\begin{tabular}{|c|c|c|c|c|c|c|}
\hline
Case&$d$& $g$& $n$&$b$&  Reference & $\dim{\mathcal{H}}$\\
\hline
\hline
1&$7$&$3$&$6$&$1$&\cite{Io1} Theorem 4.3 and \S 8&$64$\\
\hline
2&$8$&$3$&$7$&$0$&\cite{Io1} Theorem 4.3&$74$\\
\hline
3&$9$&$3$&$8$&$-1$&\cite{fa-li9} Theorem 3.3&$86$\\
\hline
4&$9$&$4$&$7$&$1$&\cite{fa-li9} Theorem 3.3&$84$\\
\hline
5&$10$&$3$&$9$&$-2$&\cite{fa-li10} Prop. 3.4&$100$\\
\hline
6&$10$&$4$&$8$&$0$&\cite{fa-li10} Prop. 3.4&$96$\\
\hline
7&$10$&$5$&$7$&$2$&\cite{fa-li10} Prop. 7.1&$94$\\
\hline
8&$11$&$3$&$10$&$-3$&\cite{be-bi} Prop. 4.2.2&$116$\\
\hline
9&$11$&$4$&$9$&$-1$&\cite{be-bi} Prop. 4.2.3, Remark 4.2.5&$110$\\
\hline
10&$11$&$5$&$8$&$1$&\cite{be-bi} Prop. 5.2.1&$106$\\
\hline
11&$11$&$6$&$7$&$3$&\cite{be-bi} Prop. 5.2.1&$104$\\
\hline

\end{tabular}
\end{center}
\end{cor}
\vskip .15 in
\begin{proof}
We will show that the hypothesis of Proposition \ref{Hilbertscheme of
fibroverP1} are satisfied for all cases in the above table .
The structure sequences of a general surface section $S$ and curve
section $C$ tensored with $L$ and $\restrict{L}{S}$ respectively, give
\begin{equation}
h^0(L) - h^1(L)= 2 + \chiS + d - g.
\end{equation}
On the other hand, $h^0(L) = n+1$ and thus
\begin{equation}
\label{h1LfibrationoverP1}
h^1(L) = n -1 - \chiS - d + g.
\end{equation}
Noticing that $\chiS = 1,$ being $S$ a conic bundle over $\Pin{1},$ a
simple check gives $h^1(L)= 0$ for all cases in the above table.

As in the proof of \cite{fug2}, Lemma 3.19, in all cases above, but
{\it 7} and {\it 11}, as  $b \le 1$, it follows that $a_1 \ge 0.$
A simple numerical check, independent of the value of $a_1\ge 0$ shows
that hypothesis ii) in Proposition \ref{Hilbertscheme of fibroverP1} is
verified in cases {\it 1} through {\it 4, 6, 9} and {\it 10}.

In cases {\it 5} and {\it 8}, \cite{fa-li1} Theorem 2.0 gives $a_1 = 1$
and thus hypothesis ii) in Proposition \ref{Hilbertscheme of
fibroverP1} is satisfied.

Because $H^1(\Pin{1},\mathcal{E}) = H^1(X,L) = 0,$ it follows that
$a_1 \ge -1.$ Therefore hypothesis ii) in Proposition
\ref{Hilbertscheme of fibroverP1} is verified in cases {\it 7} and {\it
11}.

By Proposition \ref{Hilbertscheme of fibroverP1} there exists an
irreducible component $\mathcal{H}$ of the Hilbert scheme of $X \subset
\Pin{n}$ whose dimension can now be easily computed.
\end{proof}

%
%

\begin{cor}
\label{cordelpezzo} Let $\xel$ be as in {\rm \ref{fibrationsonP1}},
with $\alpha = 3,$ i.e. a Del Pezzo fibration over $\Pin{1}$ with fiber
degree $3.$ Let $X$ be embedded by $|L|$ in $\Pin{n}$ and let the
numerical invariants of $X \subset \Pin{n}$ be as in the table below,
where $S \in |L|$ is a general surface section.
Then the  Hilbert scheme of $X \subset \Pin{n}$ has an irreducible
component, ${\mathcal
H}$,  which is smooth at the point
representing $X$ and of dimension as in the rightmost column of the
table.
\vskip .2in
\begin{center}
\begin{tabular}{|c|c|c|c|c|c|c|}
\hline
$d$& $g$& $n$&$p_g(S)$&$b$& Reference & $\dim{\mathcal{H}}$\\
\hline
\hline
$9$&$7$&$6$&$2$&$0$&\cite{fa-li9} Prop. 2.5,&$94$\\
\hline
$10$&$9$&$6$&$3$&$1$&\cite{fa-li10} Theorem 4.2,&$114$\\
\hline
$11$&$8$&$7$&$2$&$-1$&\cite{be-bi} Remark 5.4.7&$104$\\
\hline
\end{tabular}
\end{center}
\end{cor}
\begin{proof}
We will show that the hypothesis of Proposition \ref{Hilbertscheme of
fibroverP1} are satisfied for all cases in the above table.

Noticing that formula \brref{h1LfibrationoverP1} holds true in these
cases and that  $q(S) = 0,$ a simple check gives $h^1(L)= 0.$

As in the proof of \cite{fug2}, Lemma 3.19, because in all cases above
it is  $b \le 2$, it follows that $a_1 \ge 0.$
A simple numerical check now shows that hypothesis ii) in Proposition
\ref{Hilbertscheme of fibroverP1} is verified in all the cases of the
above table.

By Proposition \ref{Hilbertscheme of fibroverP1} there exists an
irreducible component $\mathcal{H}$ of the Hilbert scheme of $X \subset
\Pin{n}$ whose dimension can now be easily computed.
\end{proof}

\section{Good determinantal varieties}
\label{gooddet}
  In \cite{kmmnp} results are obtained on the unobstructedness of good
determinantal subschemes, as an application of the authors' remarkable
G-liaison theory. In particular, all good determinantal subschemes of
codimension 3 were shown to be unobstructed and the dimension of their
locus inside the Hilbert scheme was computed.
Some of the varieties considered in previous sections of this work have
explicit constructions, available in the literature, which easily show
that they are examples of good determinantal schemes. This section
addresses the relationship between our work and \cite{kmmnp}.

For the convenience of the reader we begin the section by recalling the
definition of good determinantal
subscheme and the basic notation utilized in \cite{kmmnp}.

  \begin{dfntn}[cf.Definition 3.1 in \cite{kmmnp}]
  \label{gooddetdef}
  Let $A$ be a homogenous matrix, i.e
  $A = [f_{i,j}]$ where $f_{i,j}\in \mathbb{C}[x_0, \dots x_{n+c}]$ is
  a homogeneous polynomial of degree $d_{ij}.$
  Let $I(A)$ denote the ideal of maximal minors of $A.$
  A codimension $c$ scheme, $X,$ in $\Pin{n+c}$ will be called {\it
standard  determinantal scheme} if $I_{X}=I(A)$ for some  homogeneous
  $t\times(t+c-1)$ matrix, $A$.
  $X$ will be called a {\it good determinantal scheme} if
  additionally, $A$ contains a  $(t-1)\times(t+c-1)$ submatrix (allowing
  a change of basis if necessary) whose
  ideal of maximal minors defines a subscheme of codimension $c+1$.
If $u: F\to E$ is a vector bundle
  homomorphism, over $\Pin{n+c},$ we define $I(u)=I(A)$ for any
homogeneous matrix
  representing $u.$
\end{dfntn}

Note that the matrix $A$ defines a morphism of locally free sheaves
$$u: {\oplus}_{i=1}^t {\mathcal O}_{\Pin{n+c}}(b_i) \to
{\oplus}_{j=0}^{t+c-2} {\mathcal O}_{\Pin{n+c}}(a_j)$$
where $d_{ij} = a_j - b_i,$ $b_1\ge \dots \ge b_t$ and $a_0\ge \dots
\ge a_{t+c-2}.$

The locus of good determinantal
subschemes $X$ in $\Pin{n+c}$  of  codimension
$c$ defined by a matrix $A$ as above, following \cite{kmmnp},
is denoted by
  $W(\underline{b},\underline{a})$,
where
$W(\underline{b},\underline{a})$
stands for
  $W(b_1,\dots,b_t;a_0,\dots,a_{t+c-2}),$
$W(\underline{b},\underline{a})\subset$ Hilb$^{p(t)}(\Pin{n+c}),$ where
$p(t)$
is the Hilbert polynomial of $X.$

\vspace{3mm}

The following examples are explicit constructions of some of the
varieties which appear in the Table in Corollaries \ref{corscrollp2},
  \ref{cordelpezzo} and  \ref{corquadricfibrations},
respectively.
\begin{ex*}
  \label{scrollsuP2}  Let $F = {\mathcal O}_{\Pin{6}}^{\oplus 3}$
and $E = {\mathcal O}_{\Pin{6}}(1)^{\oplus 5}$ be vector bundles on
$\Pin{6}$ and let
$u:{\mathcal O}_{\Pin{6}}^{\oplus 3} \to
{\mathcal O}_{\Pin{6}}(1)^{\oplus 5}$ be a generic
vector bundle homomorphism. Let  $I(u)=I(A)$, where $A$
is a homogeneous matrix
  representing $u$  and  where $I(A)$  denotes the ideal of
  maximal minors of $A.$ Let
  $X_1$ be the determinantal variety whose ideal $I_{X_{1}}=I(u)$.
  It can be easily seen that any such $X_{1}$ is a smooth threefold in
$\Pin{6}$
  with $\deg{X_1} =10$, $g(X_1)=6,$ Hilbert polynomial $p_1(t) =
\frac{5}{3}t^{3}+4t^{2}+\frac{10}{3}t+1$ and that it has the structure
of
  a scroll over $\Pin{2},$ see  (\cite{fa-li10}, Prop. 5.2).
  It appears in row $4$ of the Table in Corollary \ref{corscrollp2}.

It is straightforward from Definition \ref{gooddetdef} to see that
this
variety is a good determinantal subscheme of $\Pin{6}.$
\end{ex*}

\begin{ex*}
  \label{fibsuP1} Let
$E_2={\mathcal O}_{\Pin{6}}(1)^{\oplus 3}{\oplus}{\mathcal
O}_{\Pin{6}}(3)$,
$E_3={\mathcal O}_{\Pin{6}}(1)^{\oplus 3}{\oplus}{\mathcal
O}_{\Pin{6}}(2)$
  and $F={\mathcal O}_{\Pin{6}}^{\oplus 2}$
  be vector bundles over $\Pin{6}.$ Let
$u_i: F\to  E_i$, for $i=2,3$,
be  generic vector bundle homomorphisms and let $X_i$
  be the determinantal varieties whose ideals $I_{X_{i}}=I(u_{i})$.
  It can be easily seen that $\deg{X_2} =
10$ and $g(X_2)= 9$ while  $\deg{X_3} =
7$ and $g(X_3) = 3.$

Any such $X_2$ is a smooth threefold, known to be
a Del Pezzo fibration of fiber degree
$3$ over ${\Pin{1}},$ see (\cite{fa-li10}, Remark 4.3), with Hilbert
polynomial
$p_2(t)=\frac{5}{3}t^{3}+t^{2}+\frac{10}{3}t+1.$
It appears in row 2 of the Table in Corollary \ref{cordelpezzo}.

Any such $X_3$ is a smooth threefold, known to be a quadric fibration
over
${\Pin{1}},$ with Hilbert polynomial
$p_3(t)=\frac{7}{6}t^{3}+\frac{5}{2}t^{2}+\frac{7}{3}t+1.$ It appears
in the
first row of the Table in Corollary \ref{corquadricfibrations}.

In both cases it is straightforward
from Definition \ref{gooddetdef} to see that these varieties are
good determinantal subschemes of $\Pin{6}.$
\end{ex*}

\begin{prop} Let $X \subset \Pin{6}$ be a threefold of degree $d,$
genus $g,$ with Hilbert polynomial $p(t),$ as in one of the cases in
the following table
\vskip .2in
\begin{center}
\begin{tabular}{|c|c|c|c|c|c|c|}
\hline
$d$& $g$& $p(t)$ &Geometric Structure& Reference\\
\hline
\hline
& & & & \\
$7$&$3$&$\frac{7}{6}t^{3}+\frac{5}{2}t^{2}+\frac{7}{3}t+1$&
Hyperquadric fibration &Corollary \ref{corquadricfibrations},\\ & &
&over $\Pin{1}$& \cite{Io1} Theorem 4.3 and \S 8 \\
\hline
& & & & \\
$10$&$6$&$\frac{5}{3}t^{3}+4t^{2}+\frac{10}{3}t+1$&Scroll over
$\Pin{2}$&Corollary \ref{corscrollp2},\\
& & & &  \cite{fa-li10} Remark 5.3, \\
\hline
& & & & \\
$10$&$9$&$\frac{5}{3}t^{3}+t^{2}+\frac{10}{3}t+1$&DelPezzo fibration
&Corollary \ref{cordelpezzo},\\
& & &over $\Pin{1}$ of fiber degree $3$ & \cite{fa-li10} Theorem 4.2\\
\hline
\end{tabular}
\end{center}
Then:
\begin{itemize}
\item[i)] $Hilb^{p(t)}(\Pin{6})$ has an irreducible component
$\mathcal{H}$ of dimension ,respectively, $\dim{\mathcal{H}} = 64,
72,114$  and $X$ corresponds to a smooth point of $\mathcal{H};$
\item[ii)] There exists an open subset $U \in \mathcal{H}$ whose points
correspond to good determinantal subschemes given by  morphisms of
vector bundles $$u: {\oplus}_{i=1}^t {\mathcal O}_{\Pin{6}}(b_i) \to
{\oplus}_{j=0}^{t+1} {\mathcal O}_{\Pin{6}}(a_j)$$
with $b_1\ge \dots \ge b_t$ and $a_0\ge \dots \ge a_{t+1};$
\item[iii)] Let $W(\underline{b}, \underline{a}) \subset
Hilb^{p(t)}(\Pin{6})$ be the locus of good determinantal subschemes
where $\underline{a} = \{a_0, \dots,a_{t+1}\}$ and $\underline{b} =
\{b_0, \dots,b_t\}$ are as in $\rm{ii)}.$ Then
$\mathcal{H}$ is the closure of $W(\underline{b}, \underline{a})$ in
$Hilb^{p(t)}(\Pin{6})$
\end{itemize}
\end{prop}

\begin{proof}
Statement i) follows from Corollaries \ref{corscrollp2},
\ref{corquadricfibrations} and \ref{cordelpezzo}.
Statement ii) follows from Examples \ref{scrollsuP2} and \ref{fibsuP1}.
Let $X_i, i=1,2,3$ be as in Examples \ref{scrollsuP2} and
\ref{fibsuP1}. One can easily check that all of these varieties satisfy
the hypotheses of
\cite[Corollary 10.15]{kmmnp}
and thus they are unobstructed. Let $W_i$ denote the locus of good
determinantal subschemes with Hilbert polynomial $p_i(t)$ to which the
$X_i$ belong. The same corollary gives, respectively, $\dim{W_i} =
72,64,114$ for $i = 1,2,3.$

  Hence ${\mathcal H}$ and the locus $W_i$ of good
determinantal subschemes in $\Pin{6}$ have the same dimension. As each
$X_i$ is a smooth point of $\mathcal{H},$ it must be
${\mathcal H}=\overline{W}$, where $\overline{W}$ denotes the closure of
$W$ in $Hilb^{p(t)}(\Pin{6})$. Note that this is not true in general,
see
\cite[Example 10.5 (4)]{kmmnp}.

\end{proof}


\begin{thebibliography}{10}

\bibitem{ADS}
H.~Abo, W.~Decker, and N.~Sasakura.
\newblock An elliptic conic bundle in $\mathbb{P}^{4}$ arising from a
stable
   rank-3 vector bundle.
\newblock {\em Math. Z.}, 229(4):725--748, 1997.

\bibitem{BSS1}
M.~Beltrametti, M.~Schneider, and A.~J. Sommese.
\newblock Threefolds of degree $9$ and $10$ in $\mathbf{P}^5$.
\newblock {\em Math. Ann.}, (288):613--644, 1990.

\bibitem{BSS2}
M.~Beltrametti, M.~Schneider, and A.~J. Sommese.
\newblock Threefolds of degree $11$ in $\mathbf{P}^5$.
\newblock In {\em Complex Projective Geometry, Bergen - Trieste},
number 179 in
   London Math. Soc. Lecture Notes, pages 59--80, 1992.

\bibitem{BESO}
M.~Beltrametti and A.~J. Sommese.
\newblock {\em The Adjunction Theory of Complex Projective Varieties},
   volume~16 of {\em Expositions in Mathematics}.
\newblock De Gruyter, 1995.

\bibitem{ber}
M.~Bertolini.
\newblock Threefolds in $\mathbb{P}^6$ of degree twelve.
\newblock {\em Preprint}, 2004.

\bibitem{be-bi}
G.~M. Besana and A.~Biancofiore.
\newblock Degree eleven manifolds of dimension greater or equal to
three.
\newblock {\em Forum Mathematicum, to appear}, 2004.

\bibitem{bi-fa}
A.~Biancofiore and M.~L. Fania.
\newblock On the structure of linked threefolds.
\newblock {\em Rev. Mat. Complut.}, (1):17--44, 2001.

\bibitem{boss}
R.~Braun, G.~Ottaviani, M.~Schneider, and F.-O. Schreyer.
\newblock Boundedness for nongeneral-type {$3$}-folds in {${\bf P}\sb
5$}.
\newblock In {\em Complex analysis and geometry}, Univ. Ser. Math.,
pages
   311--338. Plenum, New York, 1993.

\bibitem{chang1}
M.-C. Chang.
\newblock The number of components of {H}ilbert schemes.
\newblock {\em Internat. J. Math.}, 7(3):301--306, 1996.

\bibitem{chang2}
M.-C. Chang.
\newblock Inequidimensionality of {H}ilbert schemes.
\newblock {\em Proc. Amer. Math. Soc.}, 125(9):2521--2526, 1997.

\bibitem{ds1}
H.~D'Souza.
\newblock Threefolds whose hyperplane sections are elliptic surfaces.
\newblock {\em Pacific J. Math.}, 134(1):57--78, 1988.

\bibitem{e}
G.~Ellingsrud.
\newblock Sur le sch\'ema de {H}ilbert des vari\'et\'es de codimension
2 dans
   $\mathbb{P}^e$ a c\^{o}ne de {C}ohen-{M}acaulay.
\newblock {\em Ann. scient. \'Ec. Norm. Sup. $4^{e}$ serie},
(8):423--432,
   1975.

\bibitem{fa-li1}
M.~L. Fania and E.~L. Livorni.
\newblock Polarized manifolds of dimension $\ge 3$, {$\Delta$}-genus
three, dim
   {B}s{$|L|$} $\le 0$ and degree $\ge 2${$ \Delta$}$ -1$.
\newblock {\em Saitama Math. Journal}, (11), 1993.

\bibitem{fa-li9}
M.~L. Fania and E.~L. Livorni.
\newblock Degree nine manifolds of dimension $\ge 3$.
\newblock {\em Math. Nachr.}, (169):117--134, 1994.

\bibitem{fa-li10}
M.~L. Fania and E.~L. Livorni.
\newblock Degree ten manifolds of dimension $n$ greater than or equal
to $3$.
\newblock {\em Math. Nachr.}, (188):79--108, 1997.

\bibitem{fa-me}
M.~L. Fania and E.~Mezzetti.
\newblock On the {H}ilbert scheme of {P}alatini threefolds.
\newblock {\em Advances in Geometry}, (2):371--389, 2002.

\bibitem{fug2}
T.~Fujita.
\newblock Classification of polarized manifolds of sectional genus two.
\newblock In {\em Algebraic Geometry and Commutative Algebra in honor of
   Masayoshi Nagata}, pages 73--98. Kinokuniya, 1987.

\bibitem{groth}
A.~Grothendieck.
\newblock {\em Techniques de construction et th\'eor\`emes d'existence
en
   g\'eom\'etrie alg\'ebrique. IV: les sch\'emas de {H}ilbert}.
\newblock Number 221 in Seminaire Bourbaki. 1960.

\bibitem{H}
R.~Hartshorne.
\newblock {\em Algebraic Geometry}.
\newblock Number~52 in GTM. Springer Verlag, New York - Heidelberg -
Berlin,
   1977.

\bibitem{Io1}
P.~Ionescu.
\newblock Embedded projective varieties of small invariants.
\newblock In {\em Proceedings of the Week of Algebraic Geometry
Bucharest
   1982}, number 1056 in Lecture Notes in Mathematics, pages 142--186.
Springer
   Verlag, 1984.

\bibitem{Io3}
P.~Ionescu.
\newblock Generalized adjunction and applications.
\newblock {\em Math. Proc. Camb. Phil. Soc.}, (99):457--472, 1986.

\bibitem{Io2}
P.~Ionescu.
\newblock Embedded projective varieties of small invariants {III}.
\newblock In {\em Algebraic Geometry L'Aquila 1988}, number 1417 in
Lecture
   Notes in Mathematics, pages 138--154. Springer Verlag, 1990.

\bibitem{kmmnp}
J.~O. Kleppe, J.~C. Migliore, R.~Mir{\'o}-Roig, U.~Nagel, and
C.~Peterson.
\newblock Gorenstein liaison, complete intersection liaison invariants
and
   unobstructedness.
\newblock {\em Mem. Amer. Math. Soc.}, 154(732):viii+116, 2001.

\bibitem{kl-mr}
J.~O. Kleppe and R.~M. Mir\'o-Roig.
\newblock The dimension of the {H}ilbert scheme of {G}orenstein
codimension 3
   subschemes.
\newblock {\em J. Pure Appl. Algebra}, (127):73--82, 1998.

\bibitem{ok2}
C.~Okonek.
\newblock \"{U}ber 2-codimensionale {U}ntermannigfaltigkeiten vom
{G}rad $7$ im
   $\mathbf{P}^4$ und $\mathbf{P}^5$.
\newblock {\em Math. Z.}, (187):209--219, 1984.

\bibitem{ok8}
C.~Okonek.
\newblock Fl\"achen vom {G}rad $8$ im $\mathbf {P}^4$.
\newblock {\em Math. Z.}, (191):207--223, 1986.

\bibitem{ok3}
C.~Okonek.
\newblock Notes on varieties of codimension $3$ in $\mathbf{P}^n$.
\newblock {\em Manuscripta Math.}, 84:421--442, 1994.

\bibitem{OSS}
C.~Okonek, M.~Schneider, and H.~Spindler.
\newblock {\em Vector Bundles on Complex Projective Spaces}.
\newblock Number~3 in Progress in Mathematics. Birkh\"auser, Boston -
Basel -
   Stuttgart, 1980.

\bibitem{ot1}
G.~Ottaviani.
\newblock On $3$-folds in $\mathbb{P}^{5}$ which are scrolls.
\newblock {\em Annali della Scuola Normale di Pisa Scienze Fisiche e
   Matematiche}, IV - XIX(3):451--471, 1992.

\bibitem{SO1}
A.~J. Sommese.
\newblock On the minimality of hyperplane sections of projective
threefolds.
\newblock {\em J. Reine Angew. Math.}, (329):16--41, 1981.

\bibitem{SD}
H.~Swinnerton~Dyer.
\newblock An enumeration of all varieties of degree 4.
\newblock {\em American Journal of Math}, (95):403--418, 1973.

\bibitem{We}
XXX.
\newblock Correspondence.
\newblock {\em American Journal of Mathematics}, (79):951--952, 1957.
\newblock Now in A. Weil, Oeuvres Sc. II, 555-556.

\end{thebibliography}
\end{document}